\documentclass[a4, 12pt]{amsart}
\usepackage[mathscr]{eucal}
\usepackage{amssymb}
\usepackage{latexsym}
\usepackage{amsthm}
\theoremstyle{plain}
\newtheorem{theorem}{Theorem}[section]

\newtheorem{lemma}{Lemma}[section]

\makeatletter
\@addtoreset{equation}{section}

\title[A gap theorem of self-shrinkers]{A gap theorem of self-shrinkers}
\author{Qing-Ming Cheng and  Guoxin Wei}
\address{Qing-Ming Cheng \\ Department of Applied Mathematics, Faculty of Sciences ,
Fukuoka  University, 814-0180, Fukuoka,  Japan, cheng@fukuoka-u.ac.jp}
\address{Guoxin Wei \\  School of Mathematical Sciences, South China Normal University,
510631, Guangzhou,  China, weiguoxin@tsinghua.org.cn}

\begin{document}
\maketitle

\begin{abstract}
\noindent In this paper,  we study complete self-shrinkers in Euclidean space  and prove
 that an $n$-dimensional complete self-shrinker with polynomial volume growth in Euclidean space $\mathbb{R}^{n+1}$
 is  isometric to either  $\mathbb{R}^{n}$,  $S^{n}(\sqrt{n})$, or  $\mathbb{R}^{n-m}\times  S^m (\sqrt{m})$, $1\leq m\leq n-1$,  if the squared norm $S$ of the second fundamental form is constant and satisfies $S<\frac{10}{7}$.

\end{abstract}



\footnotetext{The first author was partially  supported by JSPS Grant-in-Aid for Scientific Research (B): No. 24340013. The second author was partly supported by grant No. 11001087 of NSFC.}

\section {Introduction}
\noindent

Let $X: M\rightarrow \mathbb{R}^{n+1}$ be a smooth $n$-dimensional immersed hypersurface in the $(n+1)$-dimensional
Euclidean space $\mathbb{R}^{n+1}$. The immersed hypersurface $M$ is called a self-shrinker if it satisfies the quasilinear elliptic system:
\begin{equation*}
\mathbf{H}=-X^{N},
\end{equation*}
where  $\mathbf{H}$ denotes the mean curvature vector of $M$, $X^{N}$ denotes the orthogonal
projection of $X$ onto the normal bundle of $M$.

It it known that
self-shrinkers play an important role in the study of the mean curvature flow because
they describe all possible blow up at a given singularity of a mean curvature flow.

For $n=1$,  Abresch and Langer \cite{[AL]}  classified all smooth
closed self-shrinker curves in $\mathbb{R}^2$ and showed that the  round circle
is the only embedded self-shrinkers.
For $n\geq 2$,
Huisken \cite{[H2]}  studied compact self-shrinkers.
He proved that if  $M$ is an $n$-dimensional compact self-shrinker
with non-negative mean curvature $H$ in $\mathbb{R}^{n+1}$, then $X(M)=S^n(\sqrt{n})$.
We should notice that the condition of non-negative mean curvature is essential.
In fact, let  $\Delta$ and $\nabla$ denote the Laplacian and the gradient
operator on the self-shrinker, respectively and $\langle
 \cdot,\cdot\rangle$ denotes the standard inner product of
$\mathbb{R}^{n+1}$.
Because
$$
\Delta H-\langle X, \nabla H \rangle +SH-H=0,
$$
we obtain  $H>0$ from the maximum principle if the mean curvature is non-negative.
Furthermore,  Angenent \cite{[A]}  has constructed compact embedded self-shrinker torus
$S^1\times S^{n-1}$ in $\mathbb{R}^{n+1}$.

Huisken \cite{[H3]} and  Colding and Minicozzi \cite{[CM]} have studied complete and non-compact  self-shrinkers in  $\mathbb{R}^{n+1}$.
They have proved that  if $M$ is  an $n$-dimensional  complete  embedded self-shrinker
in $\mathbb{R}^{n+1}$ with $H\geq 0$ and  with polynomial volume growth,
then $M$ is  isometric to either the  hyperplane  $\mathbb{R}^{n}$,  the  round sphere $S^{n}(\sqrt{n})$, or a cylinder $S^m (\sqrt{m})\times \mathbb{R}^{n-m}$, $1\leq m\leq n-1$.

Without the condition $H\geq 0$,  Le and Sesum \cite{[LS]} proved
that if $M$ is an $n$-dimensional complete embedded self-shrinker  with polynomial volume growth and $S<1$
in Euclidean space $\mathbb{R}^{n+1}$, then $S=0$
and $M$ is isometric to the  hyperplane $\mathbb{R}^{n}$, where $S$ denotes the squared norm of the second fundamental form. Furthermore, Cao and  Li \cite{[CL]}  have studied the general case. They have proved
that if $M$ is an $n$-dimensional complete self-shrinker  with polynomial volume growth and $S\leq1$ in Euclidean space $\mathbb{R}^{n+1}$, then $M$ is isometric to either the hyperplane $\mathbb{R}^{n}$,  the  round sphere $S^{n}(\sqrt{n})$,  or a cylinder $S^m (\sqrt{m})\times \mathbb{R}^{n-m}$, $1\leq m\leq n-1$.

Recently,  Ding and Xin \cite{[DX1]} have studied the second gap on the squared norm of the second fundamental form and they have
proved that if $M$ is  an $n$-dimensional complete self-shrinker with polynomial volume growth in Euclidean space $\mathbb{R}^{n+1}$, there exists a positive number $\delta=0.022$ such that if $1\leq S\leq 1+0.022$, then $S=1$.

Motivated by the above results of Le and Sesum, Cao and Li, Ding and Xin, we consider the second gap for the squared norm of the second fundamental form and prove the following classification theorem for self-shrinkers:

\begin{theorem}
Let $M$ be an $n$-dimensional complete self-shrinker  with polynomial volume growth in $\mathbb{R}^{n+1}$.
If the squared norm $S$ of the second fundamental form is constant and satisfies
$$S\leq 1+\frac{3}{7},$$
 then $M$ is isometric to one of the following:

$(1)$ the hyperplane $\mathbb{R}^{n}$,

$(2)$ a cylinder $\mathbb{R}^{n-m}\times S^m (\sqrt{m})$, for $1\leq m\leq n-1$,

$(3)$ the round sphere $S^{n}(\sqrt{n})$.

\end{theorem}

\vskip 5mm

\section{Preliminaries}
\noindent
In this section, we give some notations and formulas.  Let $X: M\rightarrow \mathbb{R}^{n+1}$ be an $n$-dimensional self-shrinker in $\mathbb{R}^{n+1}$. Let $\{e_1,\cdots,e_{n},e_{n+1}\}$ be a local orthonormal basis along $M$ with dual
coframe $\{\omega_1,\cdots,\omega_{n},\omega_{n+1}\}$, such that $\{e_1,\cdots,e_{n}\}$ is a local
orthonormal basis of $M$ and $e_{n+1}$ is normal to $M$. Then we have
$$
\omega_{n+1}=0,\ \ \omega_{n+1i}=\sum_{j=1}^nh_{ij}\omega_j,\ \ h_{ij}=h_{ji},
$$
where $h_{ij}$ denotes the component of the second fundamental form of $M$.
$\mathbf{H}=\sum_{j=1}^nh_{jj}e_{n+1}$ is the mean curvature vector field, $H=|\mathbf{H}|=\sum_{j=1}^nh_{jj}$ is the mean curvature and $II=\sum_{i,j}h_{ij}\omega_i\otimes\omega_je_{n+1}$ is the second fundamental form of $M$.
The Gauss equations and Codazzi equations are given by
\begin{equation}\label{eq:12-6-4}
R_{ijkl}=h_{ik}h_{jl}-h_{il}h_{jk},
\end{equation}
\begin{equation}\label{eq:12-6-5}
h_{ijk}=h_{ikj},
\end{equation}
where $R_{ijkl}$ is the component of curvature tensor, the covariant derivative of $h_{ij}$ is defined by
$$
\sum_{k=1}^nh_{ijk}\omega_k=dh_{ij}+\sum_{k=1}^nh_{kj}\omega_{ki}+\sum_{k=1}^nh_{ik}\omega_{kj}.
$$
Let
$$F_i=\nabla_iF,\ F_{ij}=\nabla_j\nabla_iF,\ h_{ijk}=\nabla_kh_{ij}, \ {\rm and}\ h_{ijkl}=\nabla_l\nabla_kh_{ij},$$
 where
$\nabla_j $ is the covariant differentiation operator, we have
\begin{equation}\label{eq:12-6-6}
h_{ijkl}-h_{ijlk}=\sum_{m=1}^nh_{im}R_{mjkl}+\sum_{m=1}^nh_{mj}R_{mikl}.
\end{equation}
The following elliptic operator $\mathcal{L}$ is introduced by Colding and Minicozzi in \cite{[CM]}:
\begin{equation}\label{eq:12-6-1}
\mathcal{L}f=\Delta f-<X,\nabla f>,
\end{equation}
where $\Delta$ and $\nabla$ denote the Laplacian and the gradient operator on the self-shrinker, respectively and $<\cdot,\cdot>$ denotes the standard inner product of $\mathbb{R}^{n+1}$. By a direct calculation, we have
\begin{equation}\label{eq:12-6-2}
 \mathcal{L}h_{ij}=(1-S)h_{ij},\ \  \mathcal{L}H=H(1-S),\ \ \mathcal{L}X_i=-X_i,\ \ \mathcal{L}|X|^2=2(n-|X|^2),
\end{equation}
\begin{equation}\label{eq:12-6-3}
\frac{1}{2}\mathcal{L}S=\sum_{i,j,k}h_{ijk}^2+S(1-S).
\end{equation}
If $S$ is constant, then we obtain from \eqref{eq:12-6-1} and \eqref{eq:12-6-3}
\begin{equation}\label{eq:010}
\sum_{i,j,k}h_{ijk}^2=S(S-1),
\end{equation}
hence one has either
\begin{equation}
S=0,\ \ {\rm or}\ \ S=1,\ \ {\rm or}\ \ S>1.
\end{equation}
We can choose a local field of
orthonormal frames on $M^n$ such that, at the point that we consider,
$$
h_{ij}=\left \{\aligned  \lambda_i,\ \  \quad & \text {if} \quad i=j,\\
		    0,\ \  \quad  & \text {if} \quad i \neq j.\endaligned \right.
$$
then
$$
S=\sum_{i,j}h_{ij}^2=\sum_i\lambda_i^2,
$$
where $\lambda_i$ is called the principal curvature of $M$.
From \eqref{eq:12-6-4} and  \eqref{eq:12-6-6}, we get
\begin{equation}\label{eq:12-6-7}
h_{ijij}-h_{jiji}=(\lambda_i-\lambda_j)\lambda_i\lambda_j.
\end{equation}
By a direct calculation, we obtain
\begin{equation}\label{eq:01}
\sum_{i,j,k,l}h_{ijkl}^2=S(S-1)(S-2)+3(A-2B),
\end{equation}
where $A=\sum_{i,j,k}\lambda_i^2h_{ijk}^2,\ \ B=\sum_{i,j,k}\lambda_i\lambda_jh_{ijk}^2$.

We define two functions $f_3$ and $f_4$ as follows:
$$
f_3=\sum_{i,j,k} h_{ij}h_{jk}h_{ki}=\sum_{j=1}^n\lambda_j^3,\ \  \ \ f_4=\sum_{i,j,k,l} h_{ij}h_{jk}h_{kl}h_{li}=\sum_{j=1}^n\lambda_j^4,
$$
then we have

\begin{lemma}
  Let $M$ be an $n$-dimensional complete self-shrinker without boundary and with polynomial volume growth in
  $\mathbb{R}^{n+1}$. Then
  \begin{equation}\label{eq:12-6-8}
  \mathcal{L}f_3=3(1-S)f_3+6\sum_{i,j,k} \lambda_ih_{ijk}^2,
 \end{equation}
\begin{equation}\label{eq:12-6-9}
  \mathcal{L}f_4=4(1-S)f_4+4(2A+B).
 \end{equation}
 \end{lemma}

\vskip 3pt\noindent {\it Proof}.
By the definition of $f_3$, $f_4$, $\mathcal{L}f_3$ and $\mathcal{L}f_4$, we have the
following calculations:
$$
f_{3m}=3\sum_{i,j,k} h_{ijm}h_{jk}h_{ki},
$$
$$
 f_{3mm}=3\sum_{i,j,k} h_{jk}h_{ki}h_{ijmm}+3\sum_{i,j,k}h_{ijm}h_{jkm}h_{ki}+3\sum_{i,j,k}h_{ijm}h_{jk}h_{kim},
$$

$$
\Delta f_3=\sum_mf_{3mm}=3\sum_{i,j,k} h_{jk}h_{ki}\Delta h_{ij}+6\sum_{i,j,m} \lambda_i h_{ijm}^2,
$$

$$
<X,\nabla f_3>=3\sum_{i,j,k} h_{jk}h_{ki}<X,\nabla h_{ij}>,
$$

\begin{equation*}
\aligned
 \mathcal{L}f_3&=\Delta f_3-<X, \nabla f_3>\\
&=3\sum_{i,j,k} h_{jk}h_{ki}\mathcal{L}h_{ij}+6\sum_{i,j,m} \lambda_ih_{ijm}^2\\
&=3(1-S)f_3+6\sum_{i,j,k} \lambda_ih_{ijk}^2, \endaligned
\end{equation*}
and
$$
f_{4m}=4\sum_{i,j,k,l} h_{ijm}h_{jk}h_{kl}h_{li},
$$

\begin{equation*}
\aligned
f_{4mm}&=4\sum_{i,j,k,l} h_{ijmm}h_{jk}h_{kl}h_{li}+4\sum_{i,j,k,l} h_{ijm}h_{jkm}h_{kl}h_{li}\\
&\ \ +4\sum_{i,j,k,l} h_{ijm}h_{jk}h_{klm}h_{li}+4\sum_{i,j,k,l} h_{ijm}h_{jk}h_{kl}h_{lim},
\endaligned
\end{equation*}

$$
\Delta f_4=\sum_mf_{4mm}=4\sum_{i,j,k,l} h_{jk}h_{kl}h_{li}\Delta h_{ij}+4\sum_{i,j,m} \lambda_i^2h_{ijm}^2+4\sum_{i,j,k} \lambda_i\lambda_jh_{ijk}^2
+4\sum_{i,j,m} \lambda_j^2h_{ijm}^2,
$$
$$
<X, \nabla f_4>=4\sum_{i,j,k,l} h_{jk}h_{kl}h_{li}<X, \nabla h_{ij}>,
$$

\begin{equation*}
\aligned
 \mathcal{L}f_4&=\Delta f_4-<X, \nabla f_4>\\
&=4\sum_{i,j,k,l} h_{jk}h_{kl}h_{li}\mathcal{L}h_{ij}+8\sum_{i,j,m} \lambda_i^2h_{ijm}^2+4\sum_{i,j,m} \lambda_i\lambda_j h_{ijm}^2\\
&=4(1-S)f_4+4(2A+B).
\endaligned
\end{equation*}

 $$\eqno{\Box}$$

\vskip 5mm

 \section{Some estimates}
 \noindent
In this section, we will give some estimates which are needed to prove our theorem.
From now on, we denote
$$
S-1=tS,
$$
where $t$ is a positive constant if we assume that $S$ is constant and $S>1$, then
$$(1-t)S=1, \ \ \ \ \sum_{i,j,k}h_{ijk}^2=tS^2.$$

By a direct calculation, one obtains
 \begin{equation}\label{eq:1.6}
\begin{aligned}
\sum_{i,j,k,l}h_{ijkl}^2
&\geq\sum_ih_{iiii}^2+\frac{3}{4}\sum_{i\neq j}(h_{ijij}+h_{jiji})^2+\frac{3}{4}\sum_{i\neq j}(h_{ijij}-h_{jiji})^2\\
&=\sum_ih_{iiii}^2+\frac{3}{4}\sum_{i\neq j}(h_{ijij}+h_{jiji})^2\\
&\ \ +\frac{3}{2}\biggl[S\sum_{i}\lambda_i^4-(\sum_i\lambda_i^3)^2\biggl], \end{aligned}
\end{equation}
we next have to estimate $S\sum_{i}\lambda_i^4-(\sum_i\lambda_i^3)^2$ since we want to give  the estimate of $\sum_{i,j,k,l}h_{ijkl}^2$. Define
$$
f \equiv \sum_i\lambda_i^4-\frac 1S\left(
\sum_i\lambda_i^3\right)^2=f_4-\frac 1S\bigl(f_3\bigl)^2.
$$

Firstly, we have
\begin{lemma} There is one point $x \in M$ such that the following identity holds at the point.
\begin{equation}\label{eq:4.1}
\aligned
&tS^2\left[\frac{c}{S}\left(\sum_i\lambda_i^3\right)^2-\sum_i\lambda_i^4\right]\\
=&c\left(\sum_{i,j,k}2\lambda_ih_{ijk}^2\right)\sum_i\lambda_i^3 -
(2A+B)S+3c\sum_j\left(\sum_i\lambda_i^2h_{iij}\right)^2,\endaligned
\end{equation}
where $c$ is a real number.
\end{lemma}
\vskip 3pt\noindent {\it Proof}.  Define a function
$$
F=\frac 14S\sum_i\lambda_i^4-\frac16c\left(\sum_i\lambda_i^3\right)^2
=\frac 14Sf_4-\frac16c\bigl(f_3\bigl)^2,
$$
we have from Lemma 2.1,
\begin{equation}\label{eq:101}
\aligned
\mathcal{L}F&=\mathcal{L}\bigl(\frac 14Sf_4-\frac16c\bigl(f_3\bigl)^2\bigl)\\
&=S(1-S)f_4-c\biggl((1-S)f_3^2+2\sum_{i,j,k} \lambda_ih_{ijk}^2f_3\\
&\ \ +3\sum_j(\sum_i \lambda_i^2h_{iij})^2\biggl)+(2A+B)S.
\endaligned
\end{equation}
On the other hand, we have from Stokes formula,
$$
\int_{M}\mathcal{L}Fe^{-\frac{|X|^2}{2}}dv=0,
$$
hence there is a point $x \in M$ such that
$$
S(1-S)f_4-c\biggl((1-S)f_3^2+2\sum_{i,j,k} \lambda_ih_{ijk}^2f_3+3\sum_j(\sum_i \lambda_i^2h_{iij})^2\biggl)+(2A+B)S
=0$$
at the point because of the continuity of the function.
$$\eqno{\Box}$$

Secondly, we have
\begin{lemma}
$$
f=f_4-\frac{f_3^2}{S}\geq \frac{(\lambda_1-\lambda_2)^2}{\lambda_1^2+\lambda_2^2}(\lambda_1\lambda_2)^2,
$$
where $\lambda_1=\max\limits_i\{\lambda_i\} $,  $\lambda_2=\min\limits_i\{\lambda_i\}$.
\end{lemma}
\vskip 3pt\noindent {\it Proof}. Since
$$
Sf_4-f_3^2=\frac{1}{S}\sum_i(\lambda_i^2S-f_3\lambda_i)^2,
$$
then
$$
\aligned Sf_4-f_3^2 &\ge \frac{1}{S}\left( \lambda_1^2S-f_3\lambda_1\right)^2
+\frac{1}{S}\left(\lambda_2^2S-f_3\lambda_2\right)^2\\
& =S\lambda_1^4+S\lambda_2^4+\frac{f_3^2(\lambda_1^2+\lambda_2^2)}{S}-2(\lambda_1^3+\lambda_2^3)f_3\\
&\geq S(\lambda_1^4+\lambda_2^4)-\frac{\lambda_1^2
+\lambda_2^2}{S}\frac{(\lambda_1^3+\lambda_2^3)^2S^2}{(\lambda_1^2+\lambda_2^2)^2}\\
&=\frac{S}{\lambda_1^2+\lambda_2^2}\lambda_1^2\lambda_2^2(\lambda_1-\lambda_2)^2.
\endaligned
$$
$$\eqno{\Box}$$

 Thirdly, one has
\begin{lemma}
\begin{equation}\label{eq:011} A-B \le \frac 13(\lambda_1-\lambda_2)^2tS^2(1-\alpha),
\end{equation}
where $\alpha=\dfrac{\sum_ih_{iii}^2}{\sum_{i,j,k}h_{ijk}^2}=\dfrac{\sum_ih_{iii}^2 }{tS^2}$. \end{lemma}
\vskip 3pt\noindent {\it Proof}. By means of symmetry, we have
$$
\aligned
A-B&=\sum_{i,j,k}(\lambda_i^2-\lambda_i\lambda_j)h_{ijk}^2\\
   &=\frac 13\sum_{i,j,k}(\lambda_i^2+\lambda_j^2+\lambda_k^2 -
   \lambda_i\lambda_j-\lambda_j\lambda_k-\lambda_k\lambda_i)h_{ijk}^2\\
   &=\frac 13\sum_{i,j}3(\lambda_i-\lambda_j)^2h_{iij}^2\\
&\ \ +\frac 13
   \sum_{i\neq j\neq k \neq i}(\lambda_i^2+ \lambda_j^2 +\lambda_k^2
   -\lambda_i\lambda_j -\lambda_j\lambda_k -\lambda_k\lambda_i
   )h_{ijk}^2. \endaligned
$$
Without loss of generality, we can assume that $ \lambda_i \le \lambda_j
\le \lambda_k $ and consider
$$
z=\lambda_i^2+\lambda_j^2 +\lambda_k^2 -\lambda_i\lambda_j -
\lambda_j\lambda_k -\lambda_k\lambda_i
$$
as a function of $\lambda_j$, which takes its maximum at one of the
boundary points $\lambda_i$ or $\lambda_k $. On the other hand,
$$
z_{\lambda_j=\lambda_i}=z_{\lambda_j=\lambda_k}=(\lambda_i-\lambda_k)^2
\le (\lambda_1-\lambda_2)^2.
$$
Hence we get
$$
\aligned A-B \le &\frac 13\left[\sum_{i,j}
3(\lambda_i-\lambda_j)^2h_{iij}^2 +\sum_{i\neq j\neq k \neq
i}(\lambda_1-\lambda_2)^2h_{ijk}^2\right]\\
\le &\frac 13(\lambda_1-\lambda_2)^2\left(\sum_{i,j,k}h_{ijk}^2 -\sum_i
h_{iii}^2\right). \endaligned
$$
Combining \eqref{eq:010} and the definition of $\alpha$, we get $0\leq\alpha<1$
and \eqref{eq:011}.
$$\eqno{\Box}$$

From Lemma 3.1, one knows that the estimates of $\sum_k(\sum_i\lambda_i^2h_{iik}
)^2$ and $(\sum_{i,j,k}  h_{ijk}^2\lambda_i)^2$ are needed.

\begin{lemma}
\begin{equation}\label{eq:12-5-3} \sum_k\left(\sum_i\lambda_i^2h_{iik}\right)^2\le \frac
{1+2\alpha}3tS^2f,
\end{equation}
where $\alpha=\dfrac{\sum_ih_{iii}^2 }{tS^2}$.
\end{lemma}
\vskip 3pt\noindent {\it Proof}.
Since $S=\sum_{ij}h_{ij}^2$ is constant,
we have $ \sum_i\lambda_ih_{iik}=0 $, then
$$
\aligned \sum_k\left(\sum_i\lambda_i^2h_{iik}\right)^2&= \sum_k \left[
\sum_i(\lambda_i^2-a\lambda_i)h_{iik}\right]^2 \\
 &\le \sum_i(\lambda_i^2-a\lambda_i)^2\sum_{i,k}h_{iik}^2, \endaligned
$$
for any constant $a$. Let $a=\frac 1Sf_3=\frac 1S\sum_i\lambda_i^3$, we have
\begin{equation}\label{eq:12-5-1}
\sum_k\left(\sum_i\lambda_i^2h_{iik}\right)^2\le \left[\sum_i\lambda_i^4
-\frac 1S\left(\sum_i\lambda_i^3\right)^2\right]\sum_{i,k}h_{iik}^2= f\sum_{i,k}h_{iik}^2.
\end{equation}
Since
$$
\sum_{i,j,k}h_{ijk}^2 =\sum_ih_{iii}^2+3\sum_{i\neq j}h_{iij}^2+\sum_{
i\neq j \neq k \neq i}h_{ijk}^2,
$$
then
\begin{equation}\label{eq:12-5-2}
\sum_{i,k}h_{iik}^2\leq\frac{1}{3}\biggl(\sum_{i,j,k}h_{ijk}^2+2\sum_{i}h_{iii}^2\biggl)= \frac 13(1+2\alpha)\sum_{i,j,k}h_{ijk}^2=\frac 13(1+2\alpha)tS^2,
\end{equation}
combining \eqref{eq:12-5-1} and \eqref{eq:12-5-2}, we get \eqref{eq:12-5-3}.
$$\eqno{\Box}$$

\begin{lemma}
\begin{equation}\label{eq:4.6}
\left(\sum_{i,j,k}\lambda_ih_{ijk}^2\right)^2 \le \left[\frac 13(A+2B)
-\frac {4}3\sum_k\frac 1{S+2\lambda_k^2}\left( \sum_i
\lambda_i^2h_{iik}\right)^2\right]tS^2.
\end{equation}
\end{lemma}

\vskip 3pt\noindent {\it Proof}.
A straightforward computation gives
\begin{equation*}
\aligned & \left( \sum_{i,j,k}\lambda_ih_{ijk}^2\right)^2\\
   &= \left \{\frac 13\sum_{i,j,k}\left[
   (\lambda_i+\lambda_j+\lambda_k)h_{ijk}-(a_ih_{jk}+a_jh_{ki}+a_kh_{ij}
   )\right] h_{ijk}\right \}^2 \\
   & \le\frac 19\sum_{i,j,k}\left[ (\lambda_i +\lambda_j +\lambda_k
   )h_{ijk} -(a_ih_{jk} +a_jh_{ki} +a_kh_{ij} )\right]^2
   \sum_{i,j,k}h_{ijk}^2 \\
   & =\frac 19\left[ 3(A+2B)-12 \sum_{i,k}a_k\lambda_i^2h_{iik} +
   3\sum_k(S+2\lambda_k^2)a_k^2\right] tS^2, \endaligned
\end{equation*}
for any constant $ a_k \in \mathbb{R}$.
Let $$a_k=2\dfrac{\sum_i\lambda_i^2h_{iik}}{S+2\lambda_k^2},$$
then \eqref{eq:4.6} follows.
$$\eqno{\Box}$$

\vskip 5mm

\section{Proof of Theorem 1.1}
\noindent
In this section, we will prove the Theorem 1.1. The proof has three parts. In the first part of proof, we will show that
$S>1+\frac{1}{5}=1.2$ if $S>1$. In the second part, we will prove that $S>\frac{1}{0.802}>1.24688$ if $S>\frac{6}{5}=1.2$. In the third part, we will show that $S>1+\frac{3}{7}$ if $S>\frac{1}{0.802}$.

\vskip 2mm

\vskip 3pt\noindent {\it Proof of Theorem 1.1}.
\vskip 2mm
\noindent{\bf Part I}: Claim: $S>1+\frac{1}{5}=\frac{6}{5}$ if $S>1$.
\vskip 2mm
Letting $c=2$ and applying Lemma 3.1, we get
\begin{equation}\label{eq:4.7}
\aligned 0&=(S-1)[\frac{1}{2}Sf_4-f_3^2]+(2\sum_{i,j,k}\lambda_ih_{ijk}^2)f_3\\
          &\ \ -\frac{S}{2}(2A+B)+3\sum_j(\sum_i\lambda_i^2h_{iij})^2\\
          &\leq (S-1)[\frac{1}{2}Sf_4-f_3^2]+\frac{1}{2(S-1)}(2\sum_{i,j,k}\lambda_ih_{ijk}^2)^2\\
          &\ \ +\frac{S-1}{2}f_3^2 -\frac{S}{2}(2A+B)+3\sum_j(\sum_i\lambda_i^2h_{iij})^2\\
          &\leq\frac{S-1}{2}(Sf_4-f_3^2)+S\biggl[\frac{2}{3}(A+2B)
          -\frac{8}{9S}\sum_k(\sum_i\lambda_i^2h_{iik})^2\biggl]\\
          &\ \ -\frac{S}{2}(2A+B)+3\sum_j(\sum_i\lambda_i^2h_{iij})^2\\
          &\leq \frac{S-1}{2}(Sf_4-f_3^2)-\frac{S}{6}[2A-5B]+\frac{19}{27}(1+2\alpha)tS^2f\\
          &=-\frac{S}{6}(2A-5B)+[\frac{65}{54}+\frac{38}{2}\alpha]tS^2f,
 \endaligned
\end{equation}
at the point $x$,
then it follows that
\begin{equation}\label{eq:21}
-\frac{65}{9}tSf\leq-\frac{65(2A-5B)}{65+76\alpha}.
\end{equation}
On the other hand, we have
\begin{equation}\label{eq:22}
\frac{3}{2}Sf\leq S(S-1)(S-2)+3(A-2B).
\end{equation}
Combining \eqref{eq:21} and \eqref{eq:22}, we obtain
$$
\aligned &\ \ \ \ \frac{3}{2}[1-\frac{130}{27}t]Sf\\
      &\leq S(S-1)(S-2)+3(A-2B)-\frac{65(2A-5B)}{65+76\alpha}\\
      &=S(S-1)(S-2)+4(A-B)-\frac{65(3A-3B)}{65+76\alpha}-\frac{76\alpha(A+2B)}{65+76\alpha}\\
      &\leq S(S-1)(S-2)+[4-\frac{195}{65+76\alpha}](A-2B)\ \ \ \  ({\rm Since}\ A+2B\geq 0)\\
      &\leq S(S-1)(S-2)+[4-\frac{195}{65+76\alpha}]\frac{1-\alpha}{3}(\lambda_1-\lambda_2)^2tS^2.
 \endaligned
$$
Letting $y=65+76\alpha$, we get
$$
\aligned &\ \ \ \ (4-\frac{195}{65+76\alpha})(1-\alpha)=\frac{1}{76}(4-\frac{195}{y})(141-y)\\
         &=\frac{1}{76}(564+195-\frac{195\times 141}{y}-4y)\leq \frac{1}{76}(759-2\sqrt{4\times 195\times 141})
         \equiv 3\gamma_1,
 \endaligned
$$
where $\gamma_1=0.4198\cdots<0.42$.

Since we assume $t\leq\frac{1}{6}$, that is, $1\leq S\leq 1+\frac{1}{5}=\frac{6}{5}$, then
\begin{equation}\label{eq:23}
S(S-1)(S-2)+\frac{3\gamma_1}{3}(\lambda_1-\lambda_2)^2tS^2\geq\frac{3}{2}(1-\frac{130}{27}t)
      \frac{(\lambda_1-\lambda_2)^2}{\lambda_1^2+\lambda_2^2}(\lambda_1\lambda_2)^2S.
\end{equation}
We next consider two cases:

\noindent{\bf Case 1}: $\lambda_1(x)\lambda_2(x)\geq0$.

We see from \eqref{eq:23} that
$$
S(S-1)(S-2)\geq-\gamma_1(\lambda_1-\lambda_2)^2tS^2\geq -S\gamma_1tS^2,
$$
that is,
$$
S-2\geq-\gamma_1S,
$$
then
$$
S\geq\frac{2}{1+\gamma_1}\geq\frac{2}{1+0.42}>1.4>1.2=\frac{6}{5}.
$$

\noindent {\bf Case 2}: $\lambda_1(x)\lambda_2(x)<0$.

From \eqref{eq:23}, we obtain
$$
\aligned (S-1)(S-2)+\gamma_1StS&\geq (S-1)(S-2)+\gamma_1(\lambda_1^2+\lambda_2^2)tS\\
       &\geq 2\gamma_1\lambda_1\lambda_2tS+\frac{3}{2}(1-\frac{130}{27}t)(\lambda_1\lambda_2)^2\\
       &\geq 2\gamma_1\lambda_1\lambda_2tS+\frac{3}{2}(1-\frac{130}{27}\frac{1}{6})(\lambda_1\lambda_2)^2\\
       &\geq-\frac{4\gamma_1^2t^2S^2}{4[\frac{3}{2}\times(1-\frac{130}{27\times 6})]},
 \endaligned
$$
that is,
$$
S\geq\frac{16+27\gamma_1^2}{8+8\gamma_1+27\gamma_1^2}>1.286>1.2=\frac{6}{5}.
$$
Hence we have proved
$$
S>1+\frac{1}{5}=\frac{6}{5}.
$$
\\
\noindent{\bf Part II}: Claim: $S>\frac{1}{0.802}>1.24688$ if $S>\frac{6}{5}$.
\vskip 2mm
Letting $c=\frac{9}{5}$ and applying Lemma 3.1, we have
\begin{equation}\label{eq:24}
\aligned 0&=(S-1)[\frac{5}{9}Sf_4-f_3^2]+(2\sum_{i,j,k}\lambda_ih_{ijk}^2)f_3\\
          &\ \ -\frac{5S}{9}(2A+B)+3\sum_j(\sum_i\lambda_i^2h_{iij})^2\\
          &\leq (S-1)[\frac{5}{9}Sf_4-f_3^2]+\frac{9}{16(S-1)}(2\sum_{i,j,k} \lambda_ih_{ijk}^2)^2\\
          &\ \ +\frac{4(S-1)}{9}f_3^2 -\frac{5S}{9}(2A+B)+3\sum_j(\sum_i\lambda_i^2h_{iij})^2\\
          &\leq\frac{5tS}{9}(Sf)+\frac{9}{16tS}(2\sum_{i,j,k} \lambda_ih_{ijk}^2)^2-\frac{5}{9}S(2A+B)
          +3\sum_j(\sum_i\lambda_i^2h_{iij})^2\\
          &\leq\frac{5}{9}tS^2f+\frac{1}{9\times 16tS}(2\sum_{i,j,k} \lambda_ih_{ijk}^2)^2-\frac{5}{9}S(2A+B)\\
          &\ \ +
          \frac{5}{9}\biggl[\frac{4}{3}(A+2B)S-\frac{16}{9}\sum_k(\sum_i\lambda_i^2h_{iik})^2\biggl]
          +3\sum_j(\sum_i\lambda_i^2h_{iij})^2\\
          &\leq\frac{5}{9}tS^2f+\frac{1}{144tS}(2\sum_{i,j,k} \lambda_ih_{ijk}^2)^2+\frac{20}{27}(A+2B)S\\
          &\ \ -\frac{5}{9}(2A+B)S+\frac{163}{81\times 3}(1+2\alpha)tS^2f\\
          &=\frac{5}{9}tS^2f+\frac{1}{144tS}(2\sum_{i,j,k} \lambda_ih_{ijk}^2)^2-\frac{5(2A-5B)S}{27}\\
          &\ \ +\frac{163}{81\times 3}(1+2\alpha)tS^2f,
 \endaligned
\end{equation}
at the point $x$,
that is,
\begin{equation}\label{eq:25}
0\leq3tSf+\frac{3}{80tS^2}(2\sum_{i,j,k} \lambda_ih_{ijk}^2)^2-(2A-5B)+\frac{163}{81\times3}\times\frac{27}{5}(1+2\alpha)tSf,
\end{equation}
then
\begin{equation}\label{eq:26}
-\frac{298+326\alpha}{45}tSf\leq\frac{3}{80tS^2}(2\sum_{i,j,k} \lambda_ih_{ijk}^2)^2-(2A-5B).
\end{equation}
Since
\begin{equation}\label{eq:27}
\aligned (2\sum_{i,j,k}\lambda_ih_{ijk}^2)^2&\leq4S^2[\sum_i h_{iiii}^2+\frac{1}{4}\sum_{i\neq j}(h_{ijij}+h_{jiji})^2]\\
                                     &\leq4S^2[\sum_i h_{iiii}^2+\frac{3}{4}\sum_{i\neq j}(h_{ijij}+h_{jiji})^2]\\
                                     &\leq4S^2[\sum_{i,j,k,l} h_{ijkl}^2-\frac{3}{2}Sf]\\
                                     &=4S^2[S(S-1)(S-2)+3(A-2B)-\frac{3}{2}Sf],
 \endaligned
\end{equation}
then one obtains
\begin{equation}\label{eq:28}
\frac{3}{2}Sf+\frac{(2\sum_{i,j,k} \lambda_ih_{ijk}^2)^2}{4S^2}\leq S(S-1)(S-2)+3(A-2B).
\end{equation}
We now assume $\frac{1}{6}<t\leq 0.198$, that is, $S\leq\frac{1}{0.802}$, then we  will get a contradiction.

From \eqref{eq:26}, we have
\begin{equation}\label{eq:29}
-\frac{298}{225}Sf\leq\frac{9}{40S^2}(2\sum_{i,j,k} \lambda_ih_{ijk}^2)^2-\frac{298}{298+326\alpha}(2A-5B).
\end{equation}
Noting $A+2B\geq0$, we see from \eqref{eq:28} and \eqref{eq:29} that
\begin{equation}\label{eq:30}
\aligned
\frac{79}{450}Sf&\leq S(S-1)(S-2)+\biggl[4-\frac{3\times298}{298+326\alpha}\biggl](A-B)\\
                &\leq S(S-1)(S-2)+
                \biggl[4-\frac{3\times298}{298+326\alpha}\biggl]\frac{(\lambda_1-\lambda_2)^2}{3}tS^2(1-\alpha).
\endaligned
\end{equation}
On the other hand,
\begin{equation}\label{eq:31}
\aligned
&\ \ \ \ \frac{1}{3}\biggl(4-\frac{3\times298}{298+326\alpha}\biggl)(1-\alpha)\\
&=\frac{1}{3\times326}\biggl(4-\frac{3\times298}{Z}\biggl)(624-z)\\
&=\frac{1}{3\times326}(2496+894-4z-\frac{3\times298\times624}{Z})\\
&\leq\frac{1}{3\times326}(2496+894-2\sqrt{4\times3\times298\times624})\\
&\equiv\gamma_2=0.41146\cdots<0.4115,
\endaligned
\end{equation}
where $Z=298+326\alpha$.

From \eqref{eq:30}, we have
\begin{equation}\label{eq:32}
\aligned
0&\leq S(S-1)(S-2)+\gamma_2(\lambda_1-\lambda_2)^2S^2t-\frac{79}{450}Sf\\
 &\leq S(S-1)(S-2)+\gamma_2(\lambda_1-\lambda_2)^2S^2t\\
 &\ \ -\frac{79}{450}S\frac{(\lambda_1-\lambda_2)^2}{\lambda_1^2+\lambda_2^2}(\lambda_1\lambda_2)^2,
\endaligned
\end{equation}
then it follows that
\begin{equation}\label{eq:33}
S(S-1)(S-2)\geq(\lambda_1-\lambda_2)^2\biggl(-\gamma_2tS^2+\frac{79}{450}(\lambda_1\lambda_2)^2\biggl).
\end{equation}
We next consider two cases:

\noindent{\bf Case 1}: $\lambda_1(x)\lambda_2(x)>0$.

From \eqref{eq:33}, we have
$$
S(S-1)(S-2)\geq (\lambda_1-\lambda_2)^2(-\gamma_2tS^2)\geq-S\gamma_2tS^2
$$
that is,
$$
S-2\geq-\gamma_2S,
$$
then
$$
S\geq\frac{2}{1+\gamma_2}\geq\frac{2}{1+0.42}>\frac{1}{0.802}.
$$

\noindent {\bf Case 2}: $\lambda_1(x)\lambda_2(x)\leq0$.

From \eqref{eq:32}, we obtain
$$
\aligned
S(S-1)(S-2)+\gamma_2StS^2&\geq S(S-1)(S-2)+\gamma_2(\lambda_1^2+\lambda_2^2)S^2t\\
&\geq \frac{79}{450}S\frac{(\lambda_1-\lambda_2)^2}{\lambda_1^2+\lambda_2^2}(\lambda_1\lambda_2)^2
+\gamma_2(2\lambda_1\lambda_2)tS^2\\
&\geq\frac{79}{450}S(\lambda_1\lambda_2)^2+2\gamma_2\lambda_1\lambda_2tS^2\\
&\geq-\frac{(2\gamma_2tS^2)^2}{4\times \frac{79}{450}S}=-\frac{450}{79}\gamma_2^2t^2S^3,
 \endaligned
$$
that is,
$$
S\geq\frac{2+\frac{450}{79}\gamma_2^2}{1+\gamma_2+\frac{450}{79}\gamma_2^2}>1.247456>1.2469>\frac{1}{0.802}.
$$
It is a contradiction, hence we have proved
$$
S>\frac{1}{0.802}.
$$
\\

\noindent{\bf Part III}: Claim:  $S>\frac{10}{7}$ if $S>\frac{1}{0.802}$.
\vskip 2mm

Before we prove the above Claim, we will prove the following Lemma.
\begin{lemma}
 Let $M$ be an $n$-dimensional complete self-shrinker without boundary and with polynomial volume growth in
 $\mathbb{R}^{n+1}$. If the squared norm $S$ of the second fundamental form is constant, then for any constant
 $\delta>0$, $c_0\geq0$ and $c_1$
 satisfying
 \begin{equation}\label{eq:12-5-10}
 (\beta+t)c_0\delta=(\delta-1+\delta c_0)^2,
 \end{equation}
 and $\beta\geq0$, there exists a point $p_0\in M$ such that, at $p_0$,

 \begin{equation}\label{eq:100}
 \aligned
 &\ \ tS^2(S-2)\\
 &\geq (2-\delta t+c_1\delta)Sf-(5-2\delta+c_1\delta+\frac{\beta}{3})A+(6+\delta+2c_1\delta-\frac{2}{3}\beta)B\\
 &\ \ +\biggl[4\sqrt{\frac{2\beta}{3t}}-\frac{2}{t}-3(1+c_0)\delta\biggl]\frac{1}{S}
 \sum_k(\sum_i\lambda_i^2h_{iik})^2.
  \endaligned
 \end{equation}
 \end{lemma}
 \noindent
 {\it Proof}. From \cite{[DX1]}, we have

 \begin{equation}\label{eq:102}
 \int_M (A-2B-Sf)e^{-\frac{|X|^2}{2}}dv=0,
 \end{equation}
 then for any constant $c_1$, we have
 \begin{equation}\label{eq:103}
 \int_M c_1S(A-2B)e^{-\frac{|X|^2}{2}}dv=\int_M c_1S^2fe^{-\frac{|X|^2}{2}}dv \end{equation}
 since $S$ is constant. From \eqref{eq:101}, we have
 \begin{equation}\label{eq:104}
 \aligned
 &\ \ \int_M(1-S)(cf_3^2-Sf_4)e^{-\frac{|X|^2}{2}}dv\\
 &=\int_M[(2A+B)S-2cf_3\sum_{i,j,k}\lambda_ih_{ijk}^2-3c\sum_j(\sum_i\lambda_i^2h_{iij})^2]e^{-\frac{|X|^2}{2}}dv,
 \endaligned
 \end{equation}
 then
 \begin{equation}\label{eq:105}
 \aligned
 &\ \ \int_M (c_1S^2f-tS^2f_4+ctSf_3^2)e^{-\frac{|X|^2}{2}}dv\\
 &=\int_M\biggl\{2cf_3\sum_{i,j,k}\lambda_ih_{ijk}^2-(2A+B)S+3c\sum_j(\sum_i\lambda_i^2h_{iij})^2\\
 &\ \ +c_1S(A-2B)\biggl\}e^{-\frac{|X|^2}{2}}dv,
 \endaligned
 \end{equation}
 thus we have that there exists a point $p_0\in M$ such that, at $p_0$,
 \begin{equation}\label{eq:106}
 \aligned
 &\ \ c_1S^2f-tS^2f_4+ctSf_3^2\\
 &=2cf_3\sum_{i,j,k}\lambda_ih_{ijk}^2-(2A+B)S+3c\sum_j(\sum_i\lambda_i^2h_{iij})^2+c_1S(A-2B),
 \endaligned
 \end{equation}
 then,
 \begin{equation}\label{eq:107}
 \aligned
 &\ \ c_1S^2f-tS(Sf_4-f_3^2)\\
 &=(1-c)tSf_3^2+2cf_3\sum_{i,j,k}\lambda_ih_{ijk}^2-(2A+B)S\\
 &\ \ +3c\sum_j(\sum_i\lambda_i^2h_{iij})^2+c_1S(A-2B).
 \endaligned
 \end{equation}
 Putting $c=1+c_0$ with $c_0\geq0$, we get
 \begin{equation}\label{eq:108}
 \aligned
  (c_1S^2-tS^2)f&=-c_0tSf_3^2+2(c_0+1)f_3\sum_{i,j,k}\lambda_ih_{ijk}^2-(2A+B)S\\
 &\ \ +3(1+c_0)\sum_j(\sum_i\lambda_i^2h_{iij})^2+c_1S(A-2B).
 \endaligned
 \end{equation}
 For any positive constant $\delta>0$, we have from \eqref{eq:107},
 \begin{equation}\label{eq:109}
 \aligned
 \delta tSf&=c_1\delta Sf+c_0\delta tf_3^2-2(1+c_0)\frac{f_3}{S}\delta \sum_{i,j,k}\lambda_ih_{ijk}^2\\
           &\ \ +(2\delta-c_1\delta)A+(\delta+2c_1\delta)B\\
           &\ \ -3(1+c_0)\delta\frac{1}{S}\sum_j(\sum_i\lambda_i^2h_{iij})^2.
 \endaligned
 \end{equation}

 Putting
 \begin{equation}\label{eq:110}
 u_{ijkl}=\frac{1}{4}(h_{ijkl}+h_{jkli}+h_{klij}+h_{lijk}),
 \end{equation}
 by a direct computation, we have
 \begin{equation}\label{eq:111}
 \aligned
 \sum_{i,j,k,l}h_{ijkl}^2&\geq \sum_{i,j,k,l}u_{ijkl}^2+\frac{3}{4}\sum_{i,j}(h_{iijj}-h_{jjii})^2\\
           &=\sum_{i,j,k,l}u_{ijkl}^2+\frac{3}{4}\sum_{i,j}(\lambda_i-\lambda_j)^2\lambda_i^2\lambda_j^2\\
           &=\sum_{i,j,k,l}u_{ijkl}^2+\frac{3}{2}(Sf_4-f_3^2).
 \endaligned
 \end{equation}
 From \eqref{eq:01} and \eqref{eq:111}, we obtain
 \begin{equation}\label{eq:112}
 S(S-1)(S-2)+3(A-2B)\geq \sum_{i,j,k,l}u_{ijkl}^2+\frac{3}{2}Sf.
 \end{equation}
By a direct calculation, we can get
 \begin{equation}\label{eq:113}
 \aligned
 \sum_{i,j,k,l}u_{ijkl}^2&\geq\frac{1}{2}Sf+\frac{4}{tS^2}\sum_i\lambda_i^2(\sum_j\lambda_j^2h_{jji})^2-2A\\
           &\ \ +\frac{2f_3}{S}\sum_{i,j,k}\lambda_ih_{ijk}^2+\frac{1}{S^2}(\sum_{i,j,k}\lambda_ih_{ijk}^2)^2.
           \endaligned
 \end{equation}
 Combining \eqref{eq:112} and \eqref{eq:113}, we have
 \begin{equation}\label{eq:114}
 \aligned
 &\ \ S(S-1)(S-2)+3(A-2B)\\
 &\geq 2Sf-2A+\frac{4}{tS^2}\sum_i\lambda_i^2(\sum_j\lambda_j^2h_{jji})^2\\
 &\ \ +\frac{2f_3}{S}\sum_{i,j,k}\lambda_ih_{ijk}^2+\frac{1}{S^2}(\sum_{i,j,k}\lambda_ih_{ijk}^2)^2.
 \endaligned
 \end{equation}

 From \eqref{eq:109}, one has

 \begin{equation}\label{eq:115}
 \aligned
 &\ \ \delta tSf+\frac{4}{tS^2}\sum_i\lambda_i^2(\sum_j\lambda_j^2h_{jji})^2
       +\frac{2f_3}{S}\sum_{i,j,k}\lambda_ih_{ijk}^2+\frac{1}{S^2}(\sum_{i,j,k}\lambda_ih_{ijk}^2)^2\\
 &=c_1\delta Sf+c_0\delta tf_3^2-\frac{2[(1+c_0)\delta-1]f_3}{S} \sum_{i,j,k}\lambda_ih_{ijk}^2\\
 &\ \ +\frac{1}{S^2}(\sum_{i,j,k}\lambda_ih_{ijk}^2)^2 +(2\delta-c_1\delta)A+(\delta+2c_1\delta)B\\
 &\ \ +\sum_j\bigl(\frac{4\lambda_j^2}{tS}-3(1+c_0)\delta\bigl)\frac{1}{S}(\sum_i\lambda_i^2h_{iij})^2\\
 &\geq c_1\delta
 Sf+\biggl[t-\frac{[(1+c_0)\delta-1]^2}{c_0\delta}\biggl]\frac{1}{tS^2}(\sum_{i,j,k}\lambda_ih_{ijk}^2)^2\\
 &\ \ +(2\delta-c_1\delta)A+(\delta+2c_1\delta)B\\
 &\ \ +\frac{1}{S}\sum_j(\frac{4\lambda_j^2}{tS}-3(1+c_0)\delta)(\sum_i\lambda_i^2h_{iij})^2.
 \endaligned
 \end{equation}

 Taking $\delta$ and $c$, such that,
 \begin{equation}\label{eq:116}
 (\beta+t)c_0\delta=[(c_0+1)\delta-1]^2,
 \end{equation}
 with $\beta\geq0$, we have from Lemma 3.5
 \begin{equation}\label{eq:117}
 \aligned
 &\ \ \delta tSf+\frac{4}{tS^2}\sum_i\lambda_i^2(\sum_j\lambda_j^2h_{jji})^2
       +\frac{2f_3}{S}\sum_{i,j,k}\lambda_ih_{ijk}^2+\frac{1}{S^2}(\sum_{i,j,k}\lambda_ih_{ijk}^2)^2\\
 &\geq c_1\delta Sf-\beta\biggl[\frac{1}{3}(A+2B)-\frac{4}{3}\sum_k\frac{1}{S+2\lambda_k^2}(\sum_i\lambda_i^2h_{iik})^2
 \biggl]\\
 &\ \ +(2\delta-c_1\delta)A+(\delta+2c_1\delta)B\\
 &\ \ +\frac{1}{S}\sum_j\biggl(\frac{4\lambda_j^2}{tS}-3(1+c_0)\delta\biggl)(\sum_i\lambda_i^2h_{iij})^2\\
 &=c_1\delta Sf+(2\delta-c_1\delta-\frac{\beta}{3})A+(\delta+2c_1\delta-\frac{2\beta}{3})B\\
 &\ \ +\sum_k\biggl[\frac{4}{3}\beta\frac{1}{S+2\lambda_k^2}+\frac{4\lambda_k^2}{tS^2}-\frac{3(1+c_0)\delta}{S}\biggl]
   (\sum_i\lambda_i^2h_{iik})^2\\
 &\geq c_1\delta Sf+(2\delta-c_1\delta-\frac{\beta}{3})A+(\delta+2c_1\delta-\frac{2\beta}{3})B\\
 & \ \
 +\biggl[4\sqrt{\frac{2\beta}{3t}}-\frac{2}{t}-3(1+c_0)\delta\biggl]\frac{1}{S}\sum_k(\sum_i\lambda_i^2h_{iik})^2.
 \endaligned
 \end{equation}

 From \eqref{eq:114}, we have

 \begin{equation}\label{eq:118}
 \aligned
 &\ \ tS^2(S-2)+3(A-2B)\\
 &\geq2Sf-2A-\delta tSf+c_1\delta Sf+(2\delta-c_1\delta-\frac{\beta}{3})A+(\delta+2c_1\delta-\frac{2\beta}{3})B\\
 & \ \
 +\biggl[4\sqrt{\frac{2\beta}{3t}}-\frac{2}{t}-3(1+c_0)\delta\biggl]\frac{1}{S}\sum_k(\sum_i\lambda_i^2h_{iik})^2,
 \endaligned
 \end{equation}
 that is,
 \begin{equation}\label{eq:119}
 \aligned
 &\ \ tS^2(S-2)\\
 &\geq(2-\delta t+c_1\delta)Sf-(5-2\delta+c_1\delta+\frac{\beta}{3})A\\
 &\ \ +(6+\delta+2c_1\delta-\frac{2\beta}{3})B\\
 & \ \
 +\biggl[4\sqrt{\frac{2\beta}{3t}}-\frac{2}{t}-3(1+c_0)\delta\biggl]\frac{1}{S}\sum_k(\sum_i\lambda_i^2h_{iik})^2.
 \endaligned
 \end{equation}
 $$
 \eqno{\Box}
 $$

 Taking $6+\delta+2c_1\delta-\frac{2\beta}{3}=5-2\delta+c_1\delta+\frac{\beta}{3}$, we have from
\eqref{eq:12-5-10} that
$\beta=c_1\delta+3\delta+1,\  (\beta+t)c_0\delta=((c_0+1)\delta-1)^2$. Taking $\delta=\frac{17}{5},
\ c_0=\frac{6}{17}$ and applying Lemma 4.1, we obtain
$\beta=\frac{54}{5}-t$, $c_1=-\frac{2}{17}-\frac{5}{17}t$,
\begin{equation}\label{eq:200}
\aligned
&\ \ tS^2(S-2)\\
&\geq (\frac{8}{5}-\frac{22t}{5})Sf-(\frac{7}{5}-\frac{4t}{3})(A-B)\\
&\ \ +\frac{1}{S}\biggl[4\sqrt{\frac{2}{3}(\frac{54}{5t}-1)}-\frac{2}{t}-\frac{69}{5}\biggl]
\sum_k(\sum_i\lambda_i^2h_{iik})^2.
\endaligned
\end{equation}
Putting $g_1(t)=4\sqrt{\frac{2}{3}(\frac{54}{5t}-1)}-\frac{2}{t}-\frac{69}{5}$, we can obtain that
$$
g_1(t)<0,
$$
when $t>0.1978$. Since
\begin{equation}\label{eq:201}
(g_1(t))^{'}=\frac{2}{t^2}-\frac{36\sqrt{6}}{5\sqrt{-1+\frac{54}{5t}}t^2}<0
\end{equation}
when $1>t>0.14$, then we have
$g_1(t)\leq g_1(0.1978)<0$ when $0.1978<t\leq \frac{3}{10}$.

From Lemma 3.3 and Lemma 3.4, we have
\begin{equation}\label{eq:202}
\aligned
&\ \ tS^2(S-2)\\
&\geq (\frac{8}{5}-\frac{22t}{5})Sf-(\frac{7}{5}-\frac{4t}{3})(A-B)\\
&\ \ +\frac{1}{S}\biggl[4\sqrt{\frac{2}{3}(\frac{54}{5t}-1)}-\frac{2}{t}-\frac{69}{5}\biggl]
\sum_k(\sum_i\lambda_i^2h_{iik})^2\\
&\geq (\frac{8}{5}-\frac{22t}{5})Sf-(\frac{7}{5}-\frac{4t}{3})\frac{1-\alpha}{3}(\lambda_1-\lambda_2)^2tS^2\\
&\ \ +\biggl[4\sqrt{\frac{2}{3}(\frac{54}{5t}-1)}-\frac{2}{t}-\frac{69}{5}\biggl]\frac{1+2\alpha}{3}tSf\\
&=-(\frac{7}{15}-\frac{4t}{9})(1-\alpha)(\lambda_1-\lambda_2)^2tS^2\\
&\ \ +\biggl\{\biggl[\frac{8}{5t}+\frac{4}{3}\sqrt{\frac{2}{3}(\frac{54}{5t}-1)}-\frac{2}{3t}-9\biggl]\\
&\ \ +\biggl[\frac{8}{3}\sqrt{\frac{2}{3}(\frac{54}{5t}-1)}-\frac{4}{3t}-\frac{46}{5}\biggl]\alpha\biggl\}tSf.
\endaligned
\end{equation}

If $t>\frac{3}{10}$, the result is obvious true. If $t\leq\frac{3}{10}$, we will obtain a contradiction. In this case,
we have $0.198\leq t\leq \frac{3}{10}$. Putting
\begin{equation}\label{eq:203}
a(t)=\frac{8}{5t}+\frac{4}{3}\sqrt{\frac{2}{3}(\frac{54}{5t}-1)}-\frac{2}{3t}-9,
\end{equation}
\begin{equation}\label{eq:204}
b(t)=-\frac{8}{3}\sqrt{\frac{2}{3}(\frac{54}{5t}-1)}+\frac{4}{3t}+\frac{46}{5},
\end{equation}
then we have
\begin{equation}\label{eq:205}
\aligned
&\ \ tS^2(S-2)\\
&\geq-(\frac{7}{15}-\frac{4t}{9})(1-\alpha)(\lambda_1-\lambda_2)^2tS^2+[a(t)-b(t)\alpha]tSf.
\endaligned
\end{equation}
Since
$a^{'}(t)=-\frac{14}{15t^2}-\frac{12\sqrt{6}}{5\sqrt{-1+\frac{54}{5t}}t^2}<0$, we have
\begin{equation}\label{eq:206}
a(t)\geq a(\frac{3}{10})=-\frac{52}{9}+\frac{4}{3}\sqrt{\frac{70}{3}}\approx 0.662834>0.
\end{equation}
Since
$b^{'}(t)=-\frac{4}{3t^2}+\frac{24\sqrt{6}}{5\sqrt{-1+\frac{54}{5t}}t^2}>0$ if $t>0.14$, we have that $b(t)$ is
an increasing function of $t\in [0.198, \frac{3}{10}]$, then
\begin{equation}\label{eq:207}
b(t)\leq b(\frac{3}{10})=\frac{614}{45}-\frac{8}{3}\sqrt{\frac{70}{3}}\approx 0.763221>0.
\end{equation}
Therefore we get
\begin{equation}\label{eq:208}
\aligned
&\ \ tS^2(S-2)\\
&\geq-(\frac{7}{15}-\frac{4t}{9})(1-\alpha)(\lambda_1-\lambda_2)^2tS^2+[a(\frac{3}{10})-b(\frac{3}{10})\alpha]tSf.
\endaligned
\end{equation}
We next consider two cases:

\noindent{\bf Case 1}: $a(\frac{3}{10})-b(\frac{3}{10})\alpha\leq0$.

In this case $\frac{a(\frac{3}{10})}{b(\frac{3}{10})}\leq\alpha\leq1$. Since $\lambda_1, \lambda_2$ are the maximum and minimum of the principal curvatures at any point of $M$, we obtain, for any $j$,
$$
\lambda_j+\lambda_1\geq\lambda_2+\lambda_1,
$$
$$
(\lambda_1-\lambda_j)(\lambda_1+\lambda_j)\geq(\lambda_1-\lambda_j)(\lambda_1+\lambda_2).
$$
So we get
$$
\lambda_j^2-(\lambda_1+\lambda_2)\lambda_j\leq-\lambda_1\lambda_2,
$$
and
$$
f_4-(\lambda_1+\lambda_2)f_3\leq-\lambda_1\lambda_2S,
$$
then
\begin{equation}\label{eq:209}
Sf=Sf_4-f_3^2\leq-f_3^2+(\lambda_1+\lambda_2)Sf_3-\lambda_1\lambda_2S^2,
\end{equation}
\begin{equation}\label{eq:210}
Sf\leq\frac{(\lambda_1-\lambda_2)^2}{4}S^2.
\end{equation}
From \eqref{eq:208} and \eqref{eq:210}, we have
\begin{equation}\label{eq:211}
\aligned
&\ \ tS^2(S-2)\\
&\geq-(\frac{7}{15}-\frac{4t}{9})(1-\alpha)(\lambda_1-\lambda_2)^2tS^2+[a(\frac{3}{10})-b(\frac{3}{10})\alpha]
\frac{(\lambda_1-\lambda_2)^2}{4}tS^2.
\endaligned
\end{equation}
Since $a(\frac{3}{10})-b(\frac{3}{10})\alpha\leq0$, using
$-2\lambda_1\lambda_2\leq \lambda_1^2+\lambda_2^2\leq S$,
we see from \eqref{eq:211}
\begin{equation}\label{eq:212}
\aligned
S-2&\geq\biggl\{-2(\frac{7}{15}-\frac{4t}{9})+\frac{1}{2}a(\frac{3}{10})\\
&\ \ +\biggl[2(\frac{7}{15}-\frac{4t}{9})-\frac{1}{2}b(\frac{3}{10})\biggl]\alpha\biggl\}S.
\endaligned
\end{equation}
Since
\begin{equation}\label{eq:213}
\aligned
&\ \ 2(\frac{7}{15}-\frac{4t}{9})-\frac{1}{2}b(\frac{3}{10})\\
&\geq2(\frac{7}{15}-\frac{4}{9}\times \frac{3}{10})-\frac{1}{2}\times 0.77=\frac{2}{3}-\frac{1}{2}\times 0.77>0,
\endaligned
\end{equation}
we have from \eqref{eq:212}
\begin{equation}\label{eq:214}
S-2\geq-2(\frac{7}{15}-\frac{4t}{9})\biggl[1-\frac{a(\frac{3}{10})}{b(\frac{3}{10})}\biggl]S.
\end{equation}
On the other hand,
$$
\frac{a(\frac{3}{10})}{b(\frac{3}{10})}\approx 0.86847>0.86,
$$
then from \eqref{eq:214}, we see
\begin{equation}\label{eq:215}
\aligned
S-2&\geq-2(\frac{7}{15}-\frac{4t}{9})(1-0.86)S\\
   &\geq-2(\frac{7}{15}-\frac{4}{9}\times 0.198)\times 0.14 S\\
   &>0.1061S,
\endaligned
\end{equation}
hence
\begin{equation}\label{eq:216}
S>\frac{2}{1+0.1061}>1.8>\frac{10}{7}.
\end{equation}
This is impossible.

\noindent{\bf Case 2}: $a(\frac{3}{10})-b(\frac{3}{10})\alpha>0$.

In this case $\frac{a(\frac{3}{10})}{b(\frac{3}{10})}>\alpha\geq0$. From Lemma 3.2 and \eqref{eq:208}, we obtain
\begin{equation}\label{eq:217}
\aligned
&\ \ tS^2(S-2)\\
&\geq-(\frac{7}{15}-\frac{4t}{9})(1-\alpha)(\lambda_1-\lambda_2)^2tS^2+
\biggl[a(\frac{3}{10})-b(\frac{3}{10})\alpha\biggl]tSf\\
&\geq-(\frac{7}{15}-\frac{4t}{9})(1-\alpha)(\lambda_1-\lambda_2)^2tS^2\\
&\ \ +\biggl[a(\frac{3}{10})-b(\frac{3}{10})\alpha\biggl]
\frac{(\lambda_1-\lambda_2)^2}{\lambda_1^2+\lambda_2^2}(\lambda_1\lambda_2)^2tS.
\endaligned
\end{equation}
Putting $y=-\frac{\lambda_1\lambda_2}{S}$, we have $-\frac{1}{2}\leq y=-\frac{\lambda_1\lambda_2}{S}
\leq \frac{1}{2}\frac{\lambda_1^2+\lambda_2^2}{S}\leq\frac{1}{2}$, then we infer from \eqref{eq:217} that
\begin{equation}\label{eq:218}
\aligned
&\ \ S-2\\
&\geq-(\frac{7}{15}-\frac{4t}{9})(1-\alpha)(1+2y)S+\biggl[a(\frac{3}{10})-b(\frac{3}{10})\alpha\biggl](1+2y)(-y)^2S\\
&=\biggl\{\biggl[-(\frac{7}{15}-\frac{4t}{9})(1+2y)+a(\frac{3}{10})(1+2y)y^2\biggl]\\
&\ \ +\biggl[(\frac{7}{15}-\frac{4t}{9})(1+2y)-b(\frac{3}{10})(1+2y)y^2\biggl]\alpha\biggl\}S.
\endaligned
\end{equation}
Defining two functions $\rho(y)$ and $\varrho(y)$ by
\begin{equation}\label{eq:219}
\rho(y)=-(\frac{7}{15}-\frac{4t}{9})(1+2y)+a(\frac{3}{10})(1+2y)y^2,
\end{equation}
\begin{equation}\label{eq:220}
\varrho(y)=(\frac{7}{15}-\frac{4t}{9})(1+2y)-b(\frac{3}{10})(1+2y)y^2.
\end{equation}

Since $n>2$, we have $1+2y>0$, then
\begin{equation}\label{eq:221}
\aligned
\varrho(y)&=(\frac{7}{15}-\frac{4t}{9})(1+2y)-b(\frac{3}{10})(1+2y)y^2\\
          &=(1+2y)\bigl[\frac{7}{15}-\frac{4t}{9}-b(\frac{3}{10})y^2\bigl]\\
          &>(1+2y)\bigl[\frac{7}{15}-\frac{4}{9}\times \frac{3}{10}-0.7633\times\frac{1}{4}]\\
          &=(1+2y)\times 0.142508>0.
\endaligned
\end{equation}
\begin{equation}\label{eq:222}
\aligned
\rho(y)&=-(\frac{7}{15}-\frac{4t}{9})(1+2y)+a(\frac{3}{10})(1+2y)y^2\\
       &=(1+2y)[-\frac{7}{15}+\frac{4t}{9}+a(\frac{3}{10})y^2]\\
       &<(1+2y)[-\frac{7}{15}+\frac{4}{9}\times \frac{3}{10}+0.663\times \frac{1}{4}]\\
       &=(1+2y)\times (-0.1676)<0.
\endaligned
\end{equation}
By a direct calculation, we obtain
\begin{equation}\label{eq:223}
\aligned
\rho^{'}(y)&=-2(\frac{7}{15}-\frac{4t}{9})+a(\frac{3}{10})[2y+6y^2]\\
       &=-2[\frac{7}{15}-\frac{4t}{9}+a(\frac{3}{10})y+3a(\frac{3}{10})y^2]\\
       &<-2[\frac{7}{15}-\frac{4}{9}\times\frac{3}{10}-\frac{1}{12}a(\frac{3}{10})\\
       &<-2[\frac{1}{3}-\frac{1}{12}\times 0.663]<0,
\endaligned
\end{equation}
it follows that
\begin{equation}\label{eq:224}
\rho(y)\geq\rho(\frac{1}{2})=-2(\frac{7}{15}-\frac{4t}{9})+\frac{1}{2}a(\frac{3}{10}).
\end{equation}
From the above arguments, we have
\begin{equation}\label{eq:225}
\aligned
 S-2&\geq(\rho(y)+\varrho(y)\alpha)S\\
    &\geq(-2(\frac{7}{15}-\frac{4t}{9})+\frac{1}{2}a(\frac{3}{10}))S\\
    &>\bigl(-\frac{14}{15}S+\frac{8}{9}(S-1)+\frac{1}{2}\times 0.66\times S\bigl)\\
    &=-\frac{2}{45}S+0.33S-\frac{8}{9}>0.28S-\frac{8}{9},
\endaligned
\end{equation}
then
\begin{equation}\label{eq:226}
S>\frac{\frac{10}{9}}{1-0.28}>1.54>\frac{10}{7}.
\end{equation}
It is a contradiction.

Hence, we have $t>\frac{3}{10}$, that is, $S>\frac{10}{7}$ if $S>1$.
This completes the proof of Theorem 1.1.
$$\eqno{\Box}$$

\end {document}